%
%
%
%
%
%
%
 
%
%
%
%
\scrollmode
\magnification=\magstep1
\parskip=\smallskipamount

\advance\vsize-.25in
\advance\voffset.25in
\def\Id{\mathop{\rm Id}\nolimits}

\headline={\ifnum\pageno>1\smaller\ifodd\pageno\hfill A CARLEMAN TYPE
THEOREM FOR PROPER HOLOMORPHIC EMBEDDINGS \hfill\the\pageno \else
\the\pageno\hfill\uppercase{Gregery T. Buzzard and Franc
Forstneric}\hfill\fi\else\hss\fi}

\footline={\hss}

\def\demo#1:{\par\medskip\noindent\it{#1}. \rm}
\def\ni{\noindent}               
\def\npni{\hfill\break\noindent} 
\def\ll{\leftline}
\def\cl{\centerline}

%
%
\outer\def\beginsection#1\par{\bigskip
  \message{#1}\leftline{\bf\&#1}
  \nobreak\smallskip\vskip-\parskip\noindent}

%
%
\outer\def\proclaim#1:#2\par{\medbreak\vskip-\parskip
    \noindent{\bf#1.\enspace}{\sl#2}
  \ifdim\lastskip<\medskipamount \removelastskip\penalty55\medskip\fi}

\def\endpr{\hfill $\spadesuit$ \medskip}

%
%
%
%


%
%
%
%
\def\B{{\bf B}}
\def\C{{\bf C}}

\def\R{{\bf R}}

\def\Z{{\bf Z}}


%
%
%
%

\def\cC{{\cal C}}

%
%
%
\def\a{\alpha}

\def\g{\gamma}
\def\d{\delta}
\def\e{\epsilon}
\def\z{\zeta}

\def\l{\lambda}



%
%
%
%
\def\bar{\overline}              
\def\bs{\backslash}              

%
%
\def\dim{{\rm dim}\,}                    
\def\holo{holomorphic}                   
\def\ra{real-analytic}                   

\def\pc{polynomially convex}             
\def\ss{\subset\!\subset}                
\def\supp{{\rm supp}\,}                  

\def\phe{proper holomorphic embedding}

\def\hra{\hookrightarrow}

\def\Aut{{\rm Aut}}                         

\def\begin{\ll{}\vskip 10mm \nopagenumbers}  

\begin
\cl{\bf A CARLEMAN TYPE THEOREM}\medskip
\cl{\bf FOR PROPER HOLOMORPHIC EMBEDDINGS}
\bigskip\medskip
\cl{\bf Gregery T. Buzzard\footnote{}
{\rm \hskip-\parindent Buzzard's
research at MSRI is supported in part by NSF grant DMS-9022140.
Forstneric is
supported in part by an NSF grant and by a grant from the Ministry
of Science of the Republic of Slovenia.}
 and Franc Forstneric}

\font\small=cmr9
\font\smaller=cmr8
\font\ninecsc=cmcsc9

\midinsert
\narrower\narrower
\cl{\bf Abstract} 
\baselineskip=10pt
{\small In 1927, Carleman showed that a continuous, complex-valued
function on the real line can be approximated in the Whitney topology
by an entire function restricted to the real line.  In this paper, we
prove a similar result for proper holomorphic embeddings.  Namely, we
show that a proper $\cC^r$ embedding of the real line into $\C^n$ can be
approximated in the strong $\cC^r$ topology by a proper holomorphic
embedding of $\C$ into $\C^n$.
}
\endinsert

\baselineskip=\normalbaselineskip
\beginsection 1. Introduction

We denote the complex numbers by $\C$ and the real numbers by $\R$.
Let $\C^n$ be the complex Euclidean space of dimension
$n$ with complex coordinates $z=(z_1,\ldots,z_n)$.

To motivate our main result we recall the Carleman approximation
theorem [3], [10]:  {\it For each continuous function $\l\colon \R\to\C$
and positive continuous function $\eta\colon \R\to (0,\infty)$ there
exists an entire function $f$ on $\C$ such that $|f(t)-\l(t)|<\eta(t)$
for all $t\in \R$.} This extends the Weierstrass approximation theorem
in which the approximation takes place on compact intervals in $\R$.
If $\l$ is smooth, we can also approximate its derivatives by those
of $f$. A more general result was proved by Arakelian [1]
(see [13] for a simple proof).

Our main result is an extension of Carleman's theorem
to \it \phe s \rm of $\C$ into $\C^n$ for $n>1$:

\proclaim 1.1 Theorem:
Let $n>1$ and $r\ge 0$ be integers. Given a proper
embedding $\l\colon \R \hra \C^n$ of class $\cC^r$ and a continuous
positive function $\eta\colon \R \to (0, \infty)$, there exists a \phe\
$f\colon \C \hra \C^n$ such that
$$
        |f^{(s)}(t) - \l^{(s)}(t)| < \eta(t), \quad
        t \in \R, \quad 0 \le s \le r.
$$
If in addition $T = \{t_j\} \subset \R$ is discrete, there exists $f$
as above such that  
$$ f^{(s)}(t) = \l^{(s)}(t), \quad t \in T, \quad 0 \le s \le r. $$

Here $f^{(s)}(t)$ denotes the derivative of $f$ order $s$, and similarly
for $\l$.

We emphasize that, in general, one can not expect to extend
$\l$ to a holomorphic embedding of $\C$ into $\C^n$, even if
$\l$ is \ra. In this case $\l$ will extend holomorphically
to some open set in $\C$, but in general not to all of $\C$;
and even if $\l$ extends to $\C$, the extension need not be
a proper map into $\C^n$. So our result is the best possible
in this context.

Another motivation for Theorem 1.1 are the recent developments
on embedding Stein manifolds in $\C^n$; see the papers
[2], [6], [7], [8]. (For Stein manifolds and other topics
in several complex variables mentioned here we refer the reader
to H\"ormander [11].) In these papers it was shown that a
Stein manifold $M$ which admits a \phe\ in $\C^n$ for some $n>1$
also admits an embedding $f\colon M\hra \C^n$ whose image $f(M)$
contains a given discrete subset $E\subset \C^n$ (with prescribed
finite order jet of $f(M)$ at any point of $E$), and such that
$f(M)$ intersects the image of any entire map $g\colon \C^d\to\C^n$
of rank $d=n-\dim M$ at infinitely many points [6, Theorem 5.1].
This implies in particular that $f(M)$ cannot be mapped into any
hyperplane in $\C^n$ by an automorphism of $\C^n$. In this context
we recall that any Stein manifold $M$ embeds in $\C^n$ for
$n>(3\,\dim M+1)/2$ according to Eliashberg and Gromov [4].

In light of these results it is a natural question whether one can
in some sense prescribe the embedding $M\hra\C^n$ on a positive dimensional
submanifold $N\subset M$. Our result provides an affirmative answer
in the simplest case when $M=\C$ and $N=\R\times \{i0\} \subset \C$.
The details of our construction are considerable even in this simplest case,
and the full scope of the method remains to be seen. If $\l$ is of
class $\cC^\infty$, our method can be modified so that we approximate
to increasingly high order on complements of compact subsets of $\R$.
We shall not go into details of this. Another possible extension
is to approximate a proper smooth embedding by a \phe\ on a finite set
of real lines (or certain other real curves) in $\C$.

A seemingly more difficult question is whether a result of this
type holds when $N$ is a {\it complex\/} submanifold of positive
dimension in a complex manifold $M$. For instance,
\it does every \phe\ $\l\colon \C\hra \C^n$ for $n\ge 3$
extend to a \phe\ $f\colon \C^2 \hra\C^n$ ? \rm

One of the main tools in our construction is the following
result from [9]. The case $r=0$ was obtained earlier
in [8]. This result can also be obtained by methods in [5].

%
%
%
%
\proclaim 1.2 Proposition:
Let $K \subset \subset \C^n$ ($n \ge 2$) be a compact, polynomially
convex set, and let $C \subset \C^n$ be a smooth embedded arc of class
$\cC^\infty$ which is attached to $K$ in a single point of $K$.
Given a $\cC^\infty$ diffeomorphism
$F\colon K \cup C \to K \cup C' \subset \C^n$ such that $F$ is the
identity on $(K \cup C) \cap U$ for some open neighborhood $U$ of $K$,
and given numbers $r\ge 0$, $\e>0$, there exist a neighborhood $W$
of $K$ and an automorphism $\Phi \in \Aut\C^n$ satisfying
$$ \|\Phi - \Id\|_{\cC^r(W)} < \e, \qquad
   \|\Phi-F\|_{\cC^r(C)} < \e.  $$
(Here $\Id$ denotes the identity map.) Moreover, for each finite subset
$Z \subset K \cup C$ we can choose $\Phi$ such that it agrees with
the identity to order $r$ at each point $z \in Z \cap K$ and
$\Phi|_{C}$ agrees with $F$ to order $r$ at each point of $Z \cap C$.

The same result holds with any finite number of disjoint hairs
attached to $K$.

We explain our notation. $\Delta_\rho$ is the closed disk in $\C$
of radius $\rho$ and center $0$. $\B$ is the open unit ball in $\C^n$
with center $0$, and $R\B$ is the ball of radius $R$.
For a set $A \subset \C^n$ and $\rho>0$, let
$A + \rho\overline{\B} = \{ a+z \colon  a \in A, |z| \le \rho\}$.
We identify $\C$ and $\R$ with their images in $\C^n$
under the embedding $\z \to (\z, 0, \ldots, 0)$. For $1\le j\le n$
we denote by $\pi_j$ the coordinate projection
$\pi_j(z_1,\ldots,z_n)=z_j$.

In the proof we shall use special automorphisms of $\C^n$
of the form
$$ \Psi(z)=z+f(\pi z)v,\quad z\in\C^n,$$
where $v\in\C^n$, $\pi\colon \C^n\to \C^k$ is a linear map for
some $k<n$ with $\pi v=0$ (in most cases $k=1$), and $f$ is an
entire function on $\C^k$. An automorphism of this form is called
a {\it shear}; clearly $\Psi^{-1}(z)=z-f(\pi z)v$.

%
%
%
%
\beginsection 2. Some lemmas.

The following is standard, e.g., [12, prop. 2.15.4].

\proclaim 2.1 Lemma: 
Let $\l\colon\R \hra \C^n$ be a $\cC^\infty$ proper embedding.
Then there exists a continuous $\eta\colon \R \to (0,\infty)$
such that if $\gamma \colon\R \to \C^n$
with $|\g^{(s)}(t) - \l^{(s)}(t)| < \eta(t)$ for all $t \in \R$,
$s = 0, 1$, then $\g$ is a proper embedding.

Recall that a compact set $A\ss \C^n$ is \pc\ if for each
$z\in \C^n\bs A$ there is a \holo\ polynomial $P$ on $\C^n$
such that $|P(z)|> \max\{|P(w)|\colon w\in A\}$. We refer
the reader to [11] for properties of such sets.

\proclaim 2.2 Lemma:
Let $A \subset \C^n$ be compact and \pc\ and $\rho>0$.  Let $I \subset
\R$ be an interval whose endpoints lie in
$\C^n \bs (A \cup \Delta_\rho)$, and let $r, \e>0$.  Then
there exists an automorphism $\Psi(z) = z + g(z_1) e_2$ of $\C^n$
such that
\item{(i)}   $|\Psi(z) - z| < \e$ for $z \in \Delta_\rho$,
\item{(ii)}  $\|\Psi|_\R (t) - t\|_{\cC^r(I)}< \e$, and
\item{(iii)} $\Psi(t) \notin A$ for $t \in \overline{\R\bs I}$.
\npni If $Z \subset I$ is finite, we can choose $\Psi$
as above so that $g^{(s)}(t) = 0$ for $t \in Z$, $0 \le s \le r$.

\demo Proof:
Let $\mu_1< \mu_2$ denote the endpoints of $I$ in $\R$, and let
$\Gamma_j = \{(\mu_j, \zeta, 0, \ldots, 0)\colon \zeta \in \C\}$
for $j=1, 2$.  Let $R> \max\{|\mu_1|, |\mu_2|\}+1$ such that
$A\subset R\B$.  Consider the set  $E_j =A \cap \Gamma_j$.
Since $A$ is \pc, $E_j$ is \pc\ in $\Gamma_j$ and hence
$\Gamma_j\bs E_j$ is connected.  Since the endpoints of $I$ lie
in $\Gamma_j\bs E_j$, there exists a smooth curve
$\gamma_j \colon [0,1] \to \Gamma_j\bs E_j$ with
$\gamma_j(0) = (\mu_j, 0, \ldots, 0)$ and
$|\pi_2\gamma_j(1)|>R+1$ for $j=1,2$.

Since $A$ is compact, there exists $\d>0$
such that $\gamma_j([0,1]) + 3 \d \bar{\B} \subset \C^n \bs A$.
Let $\pi_2(z) = z_2$. Let
$K=\{x+iy\in \C \colon \mu_1-\d/2 \le x\le \mu_2+\d/2,\ |y|\le \rho+1\}$.
Define a function $h\colon K\cup [-R,R] \to \C$ by
$$
        h(t) = \cases{
  \pi_2 \gamma_{1}(1)&if $t \in [-R, \mu_{1}-2\d]$; \cr
  \pi_2 \gamma_{1}((\mu_1-\d-t)/\d)&if $t \in [\mu_1-2\d, \mu_1-\d]$; \cr
  \pi_2 \gamma_{2}((t-\mu_2-\d)/\d)& if $t \in [\mu_2+\d,\mu_2+2\d]$; \cr
  \pi_2 \gamma_{2}(1)& if $t \in [\mu_2+2\d, R]$; \cr
  0& otherwise.}
$$
Choose $\eta$, $0< \eta < \min\{\e,\d\}$. By Mergelyan's theorem
[14, p.386] there is an entire function $g$ on $\C$ such that
$|h(z)-g(z)|< \eta$ for $z\in K\cup [-R,R]$. The shear
$\Psi(z) = z + g(z_1) e_2$ then satisfies (i) and (iii).  Since $I
\subset {\rm Int}K$, Cauchy's estimates imply that it also satisfies
(ii) provided that $\eta>0$ is chosen sufficiently small.
The last condition on $g$ is a trivial addition to Mergelyan's theorem.
\endpr

\medskip
\proclaim 2.3 Lemma:
Let $\l\colon\R \hra \C^n$ be a proper, $\cC^\infty$ embedding,
$K \subset \C^n$ compact, $\e>0$, and $r\in\Z_+$.
Let $Z \subset \R$ be finite, and
suppose $\l(t) \in \C = \C \times \{0\}^{n-1}$ for each $t \in Z$.
Then there exists a shear $\Psi(z) = z + h(z_1) v$ for some $v \in \C^n$
with $\pi_1 v = 0$ such that 
\item{(i)}   $\Psi(\C) \cap \l(\R) = \l(Z)$,
\item{(ii)}  $|\Psi(z) - z| < \e$ for $z \in K$, and
\item{(iii)} $\Psi(z) = z + O(|z-\l(t)|^{r+1})$ as $z \to \l(t)$, for
all $t \in Z$.

\demo Proof:  
Let $Z=\{t_j\}_{j=1}^s$. Define a polynomial on $\C$
by $h(\z) = \Pi_{1\le j\le s}(\z - \pi_1 \l(t_j))^{r+1}$.
Consider the map $\Phi\colon \C \times \C^{n-1} \to \C^n$ given by
$$
        \Phi(z_1, \alpha_2, \ldots, \alpha_n) =
        (z_1, 0, \ldots,0) + h(z_1)(0, \alpha_2, \ldots, \alpha_n).
$$
Clearly $\Phi$ is an automorphism of $(\C\bs \l(Z))\times\C^{n-1}$.
Let $\Delta_{R,j}$ denote the closed disk of radius $R$ in $\C$ with
center $\pi_1 \l(t_j)$ for $j=1,2,\ldots,s$. Choose $R >0$ such that
the discs $\Delta_{R,j}$ for $1\le j\le s$ are pairwise disjoint.
Choose a $\rho$, $0<\rho<R$, such that $\rho^2$ is a regular value of
$\mu_j(t) = |\pi_1\l(t) - \pi_1 \l(t_j)|^2$ ($t\in\R$) for each
$j=1,2,\ldots, s$.

Let $M_\rho = \C \bs \cup_{1\le j\le s} {\rm int}\Delta_{\rho,j}$.
Let $\Phi_\rho = \Phi|_{M_\rho \times \C^{n-1}}$
and $\partial \Phi_\rho = \Phi|_{\partial M_\rho \times \C^{n-1}}$.
A simple check shows that $\Phi_\rho$ and $\partial \Phi_\rho$
are transverse to $\l(\R)$.  Hence by the transversality theorem,
there exists a set $A_\rho \subset \C^{n-1}$ of full measure such
that for each $\alpha=(\a_2,\ldots,\a_n) \in A_\rho$,
$\Phi(M_\rho, \a) =\{\Phi(z_1,\a) \colon z_1\in M_\rho\}$ and
$\l(\R)$ are transverse, hence disjoint by dimension considerations.

Let $A = \cap_{j=1}^\infty A_{1/j}$.  Then $A \subset \C^{n-1}$ has
full measure, and for each $\a \in A$ we see that
$\Phi(\C \bs \l(Z), \a)$ and $\l(\R)$ are disjoint.
Finally, choose $\a \in A$ such that $|h(z_1) \a| < \e$ for
$z_1 \in \pi_1(K)$, and let $\Psi(z) = z + h(z_1) \a$.
Then $\Psi(z_1,0,\ldots,0)=\Phi(z_1,\a_2,\ldots,\a_n)$,
and $\Psi$ satisfies the conclusions of the lemma.
\endpr

\proclaim 2.4 Lemma:
Let $\l\colon \R \hra \C^n$ be a $\cC^\infty$ embedding,
$f\colon\C \hra \C^n$ a \phe,
and $I \subset \R$ a closed interval with $f|_I = \l|_I$.  Let
$K \subset \C^n$ be compact and \pc, $a, r, \e>0$, and $T \subset \R$
discrete.  Suppose that $\l(t), f(t) \notin K$ for $t \in
\overline{\R\bs I}$.  Then there exists $\Phi \in
\Aut\C^n$ such that if $g = \Phi \circ f$, then 
\item{(i)} $|g^{(s)}(t) - \l^{(s)}(t)| < \e$ for $t \in [-a,a]$,
$0 \le s \le r$,
\item{(ii)} $g^{(s)}(t) = \l^{(s)}(t)$ for $t \in T \cap [-a,a]$,
$0\le s \le r$, and
\item{(iii)} $|\Phi(z) - z| < \e$ for $z \in K$.

\demo Proof:
We may assume that $I \subset (-a,a)$.  Let $I_1$, $I_2$ be the two
connected components of $\{\z \in I\colon f(\z) \in \C^n \bs K\}$
containing the respective endpoints of $I$, and let
$I_0 = I \bs (I_1 \cup I_2)$.  Let $A$ be the polynomial hull of $K
\cup f(I_0)$.  Then $A$ is the union of $K \cup f(I_0)$ and the
bounded connected components of $f(\C)\bs (K \cup f(I_0))$.  Note that
$f(I_1)$ and $f(I_2)$ lie in $f(\C) \bs A$ since $f(t)\notin K$
for all $t\in \bar{\R\bs I}$.

Let $L = A \cup f([-a,a])$.  Then $C = \overline{L \bs A}$ is the
union of two embedded arcs, each containing an endpoint of $f([-a,a])$.
Define $F$ on $L$ by $F(z) = z$ if $z \in A$, and
$F(z) = \l f^{-1}(z)$ if $z \in f([-a,a])$.  Then $F$ is a $\cC^\infty$
diffeomorphism of $L$ which extends as the identity map on
$(A \cup C) \cap U$ for some neighborhood $U$ of $A$. Apply
proposition 1.2 to get $\Phi \in \Aut \C^n$ such that
$|\Phi(z)-z| < \e$ for $z\in K$ and such that $g = \Phi \circ f$
satisfies (i) and (ii).
\endpr

\beginsection 3. Proof of theorem 1.1.

Choose a smooth cutoff function $\chi$ on $\R$ such that
$\chi(t)=1$ for $|t|$ small and $\supp\chi \subset (-1,1)$.
Define the constant $C=C_r > 1$ such that
$\|\chi h\|_{\cC^r} \leq C \|h\|_{\cC^r}$ for each
$h\in \cC^r(\R)$. We fix such $C$ for the entire proof.

By approximation we may assume that $\l\colon \R \hra \C^n$ in
theorem 1.1 is a proper $\cC^\infty$ embedding. Decreasing
$\eta$ if necessary we may also assume that $\eta$ satisfies
lemma 2.1 for $\l$ and $\eta(t) < 1/2$ for all $t\in \R$.

We use an inductive procedure to obtain a sequence of
\phe s $f_k\colon \C\hra\C^n$ such that
$f=\lim_{k\to\infty} f_k$ exists on $\C$ and satisfies theorem 1.1.
Each $f_k$ will be a restriction to $\C=\C\times\{0\}^{n-1}$
of a holomorphic automorphism of $\C^n$. The next map $f_{k+1}$
will be of the form $f_{k+1}=\Phi_{k+1}\circ f_k \circ \Psi_{k+1}$
for suitably chosen $\Phi_{k+1},\Psi_{k+1} \in\Aut\C^n$.

We will describe the case $k=1$ after the inductive step is given.
Recall that $\Delta_k$ is the closed disc
in $\C=\C\times\{0\}^{n-1}$ with center $0$ and radius $k$,
and $\B$ is the unit ball in $\C^n$. For the induction at step $k$,
suppose we have the following:

\item{(a)} closed balls $B_j = R_j \bar \B \subset \C^n$
with $R_j \ge \max\{j+1, R_{j-1}+1\}$, $j= 1, \ldots, k$,

\item{(b)} automorphisms $\Phi_1, \ldots, \Phi_k$ of $\C^n$ with
$|\Phi_j(z) - z|< 2^{-j}$ for $z \in B_{j-1}$, $j = 2, \ldots, k$,

\item{(c)} numbers $\e_j>0$ such that $\e_1< 2^{-1}$ and
$\e_j < \e_{j-1}/2 < 2^{-j}$ for $j = 1, \ldots, k$,

\item{(d)} automorphisms $\Psi_1, \ldots, \Psi_k$ of $\C^n$
of the form $\Psi_j(z)= z + g_j(z_1) e_2 + h_j(z_1)v_j$,
where $\pi_1(v_j) = 0$ and $|\Psi_j(z) - z|< \e_j$ for $|z| \le j$,

\item{(e)} closed intervals $I_j = [-a_j, a_j]$, $j = 1, \ldots, k$,
with $a_j > \max\{a_{j-1}+2, j+2\}$, and

\item{(f)} numbers $0<\d_j < C^{-1} \inf\{\eta(t)\colon t \in I_j\}$,
$j = 1, \ldots, k$,

\noindent
such that the automorphism
$$      f_k = \Phi_k \circ \cdots \Phi_1 \circ \Psi_1 \circ
        \cdots \circ \Psi_k \in\Aut\C^n
$$
(whose restriction to $\C=\C\times \{0\}^{n-1}$ provides an
embedding $\C\hra\C^n$) satisfies:

\item{($1_k$)} $f_k(\Delta_j + \e_k\bar\B) \subset {\rm Int}B_j$
\ for $j = 1, \ldots, k$,

\item{($2_k$)} $|f_k^{(s)}(t) - \l^{(s)}(t)| < \eta(t)$
for $t \in I_k$ and $0 \le s \le r$,

\item{($3_k$)} $|f_k^{(s)}(t) - \l^{(s)}(t)| < \d_k$  for
$t \in I_k \bs (-a_k+1, a_k-1)$ and $0 \le s \le r$,

\item{($4_k$)} $f_k^{(s)}(t)=\l^{(s)}(t)$ for $t \in T \cap I_k$,
$0 \le s \le r$,

\item{($5_k$)} $f_k(\C) \cap \l(\R) = \l(T \cap I_k)$,

\item{($6_k$)} $|\l(t)| > R_k+1$ for $|t|\ge a_k-1$,

\item{($7_k$)} $|f_k(t)| > R_k$ for $|t|\ge a_k-1$.

We will now show how to obtain these hypotheses at step $k+1$.
Let $I_k^1$ and $I_k^2$ be the two connected components of the set
$\{\z \in I_k\bs \Delta_{k+1} \colon |f_k(\z)| > R_k\}$
containing the respective endpoints of the interval $I_k$.
Let $I_k^0=I_k\bs(I_k^1 \cup I_k^2)$ be the middle interval.
By $(7_k)$ we have $I_k^0 \subset (-a_k+1,a_k-1)$.

Let $K_k$ be the polynomial hull of the set
$B_k \cup f_k(\Delta_{k+1} \cup I_k^0)$. Since $f_k(\C)$ is
a complex submanifold of $\C^n$ and $B_k$ is \pc, it is seen
easily that $K_k$ is contained in $B_k\cup f_k(\C)$, and it
is the union of $B_k \cup f_k(\Delta_{k+1} \cup I_k^0)$
and the bounded connected components of the complement
$f_k(\C)\bs (B_k \cup f_k(\Delta_{k+1} \cup I_k^0))$
(see Lemma 5.4 in [6]). Note that $(7_k)$ and (e) imply that
$f_k(\R\bs (-a_k+1, a_k-1)) \subset \C^n \bs K_k$.
 
Choose $R_{k+1} > R_k + 1$ such that
$K_k \subset (R_{k+1} -1)\B$, and let
$B_{k+1} = R_{k+1} \bar \B$.  Choose $a_{k+1}>a_k +2$
to get $(6_{k+1})$. We now want to approximate $\l$ on the larger
interval $I_{k+1} = [-a_{k+1}, a_{k+1}]$ by the image of the next
embedding $\C\hra\C^n$ (to be constructed). In order to apply
lemma 2.4 we first approximate $\l$ as follows:

\proclaim 3.1 Lemma:
There exists a proper $\cC^\infty$ embedding $\l_k\colon\R \hra \C^n$
satisfying
\item{(i)}  $\l_k = f_k$ on $[-a_k+1, a_k-1]$,
\item{(ii)} $\l_k = \l$ on $\R\bs I_k$,
\item{(iii)} $|\l_k^{(s)}(t) - \l^{(s)}(t)| < \eta(t)$ for
       $t \in I_k \bs (-a_k+1, a_k-1)$, $0 \le s \le r$,
\item{(iv)} $\l_k^{(s)}(t) = \l^{(s)}(t)$ for $t \in T$, $0\le s\le r$, and
\item{(v)}  $\l_k(t)\notin K_k$ when $|t|\ge a_k-1$.
 
\demo Proof:
We define the cutoff function $\chi_k$ on $\R$ using $\chi$,
so that $\chi_k =1$ on $[-a_k+1, a_k-1]$, $\chi_k =0$ on
$\R \bs I_k$, and $\|\chi_k h\|_{\cC^r} < C \|h\|_{\cC^r}$
as before. Let
$$
        \hat{\l}_k(t) = f_k(t) \chi_k(t) + \l(t) (1-\chi_k(t)),
        \quad t\in \R.
$$
By lemma 2.1, $(3_k)$, $(4_k)$, and choice of $\eta$ and $\d_k$,
we see that (i)--(iv) are satisfied for $\hat{\l}_k$ in place of $\l_k$.

To obtain (v) we use a transversality argument to perturb
$\hat{\l}_k$ on $I_k\bs (-a_k+1,a_k-1)$.
First note that if $|t|> a_k$, then
$|\hat{\l}_k(t)| = |\l(t)| > R_k+1$ by $(6_k)$, so
$\hat{\l}_k(t) \notin B_k$.  Also, by $(5_k)$, we see that
$\hat{\l}_k(t) \notin f_k(\C)$, so $\hat{\l}_k(t) \notin K_k$.  Next,
if $t \in T \cap (I_k \bs (-a_k+1, a_k-1))$, then by $(4_k)$, $(7_k)$,
and (e) we see that $\hat{\l}_k(t) = f_k(t) \notin K_k$.  Hence there
exists a neighborhood $V$ of $T \cap (I_k \bs (-a_k+1, a_k-1))$ such
that $\hat{\l}_k(\overline{V}) \cap K_k = \emptyset$.  

Thus we need only perturb $\hat{\l}_k$ on
$I_k \bs (V \cup (-a_k+1, a_k-1))$ to get (v).
Note that if $t \in I_k \bs (-a_k+1, a_k-1)$, then from
$(6_k)$ and $(2_k)$ we see that $|\hat{\l}_k(t)| > R_k+1/2$, so
$\hat{\l}_k(t) \notin B_k$.  Finally, a simple transversality argument
implies that we can make an arbitrarily small $\cC^\infty$
perturbation of $\hat{\l}_k$ to avoid $f_k(\C)$, and hence we get
$\l_k$ with $\l_k = \hat{\l}_k$ outside $I_k \bs (V \cup (-a_k+1,
a_k-1))$ and $\l_k$ satisfying (i)--(v).
\endpr

Now we can use lemma 2.4 to approximate $\l_k$, hence to approximate
$\l$. Set
$$
        \eqalign{ \d_{k+1} &= \min\{\eta(t)\colon t \in I_{k+1}\}/(2C), \cr
        \sigma_{k+1} &= \min\{\eta(t) - |\l_k^{(s)}(t) - \l^{(s)}(t)|
        \colon t \in I_{k+1}, 0 \le s \le r\} >0. \cr}
$$
Choose $\e>0$ so small that
$$
   \e < \min\{ 2^{-(k+1)}, \d_{k+1}, \sigma_{k+1} \}, \qquad
   f_k(\Delta_j + \e_k \bar{\B}) + \e \bar{\B} \subset {\rm Int}B_j,
   \ \ 1\le j\le k.
$$
Apply lemma 2.4 with $\l = \l_k$, $f = f_k$,
$I = [-a_k+1, a_k-1]$, $K = K_k$, $a = a_{k+1}$, $r$ and $T$ unchanged,
and $\e$ as above. This provides $\Phi_{k+1} \in \Aut\C^n$
and $G = \Phi_{k+1} \circ f_k \in \Aut\C^n$ satisfying
$$
        \cases{
        |\Phi_{k+1}(z) - z|< \e &  for $z \in K_k$, hence on $B_k$; \cr
        |G^{(s)}(t) - \l_k^{(s)}(t)| < \e & for $t \in I_{k+1}$,
        $0 \le s \le r$; \cr
             G^{(s)}(t)=\l_k^{(s)}(t) & for $t \in T \cap I_{k+1}$,
             $0 \le s \le r$. \cr}
$$
In particular, $(2_{k+1})$--$(4_{k+1})$ hold with $G$ in place of
$f_{k+1}$.  

Since $f_k(\Delta_{k+1}) \subset K_k \subset (R_{k+1}-1)\B$, we can
choose $\e_{k+1}'< \e_k$ small enough that $(1_{k+1})$ holds
with $G$ in place of $f_{k+1}$ and $\e_{k+1}'$ in place of $\e_{k+1}$,
and such that if $\psi \in \Aut\C^n$ with
$\|\psi(t) - t\|_{\cC^r(I_{k+1})} < \e_{k+1}'$,
then $(2_{k+1})$ and $(3_{k+1})$ hold with $G \circ \psi$ in
place of $f_{k+1}$.  Let $\e_{k+1} = \e_{k+1}'/2$.
Then with $G$ in place of $f_{k+1}$, we have $(1_{k+1})$--$(4_{k+1})$,
$(6_{k+1})$, and $G(-a_{k+1}), G(a_{k+1}) \notin B_{k+1}$
by $(6_{k+1})$ and $(2_{k+1})$.

\medskip
Next we want to obtain $(7_{k+1})$.  We do this using lemma 2.2 to
change the embedding so that the image of $\R \bs I_{k+1}$ misses
$B_{k+1}$ while leaving the embedding essentially unchanged on
$\Delta_{k+1} \cup I_{k+1}$.  
Apply lemma 2.2 with $A = G^{-1}(B_{k+1})$, $\rho = k+1$,
$I=I_{k+1}$, $r$ unchanged, $Z = T \cap I_{k+1}$ and $\e=\e_{k+1}/2$.
This gives a shear
$$
        \psi_{k+1}(z) = z + g_{k+1}(z_1) e_2
$$
with
$$ \cases{ |\psi_{k+1}(z) - z| < \e_{k+1}/2,\ \
                     z \in \Delta_{k+1}; \cr
      \|\psi_{k+1}|_\R(t) - t \|_{\cC^r(I_{k+1})} < \e_{k+1}/2; \cr
      g_{k+1}^{(s)}(t) = 0, \ \ t \in T \cap I_{k+1}, 0 \le s \le r; \cr
      \psi_{k+1}(t) \notin G^{-1}(B_{k+1}),
         \ \ t \in \bar{\R \bs I_{k+1}}. \cr}
$$
Let $H = G \psi_{k+1}$. Then with $H$ in place of $f_{k+1}$, we have
$(1_{k+1})$--$(4_{k+1})$, $(6_{k+1})$, and $(7_{k+1})$.

For the final correction, we use lemma 2.3 to obtain $(5_{k+1})$
while maintaining the other properties. Let $R > a_{k+1}$ be such
that $A=G^{-1}(B_{k+1}) \subset R\B$.  Let $\d>0$ be such that
$$
        \psi_{k+1} \left(
        [-R,R]\bs (-a_{k+1},a_{k+1}) + \d \bar\B \right)
        \cap A = \emptyset,
$$
and such that if $\theta \in \Aut\C^n$, with $|\theta(z) - z| <
\delta$ on $R \overline{\B}$, then
$$
	\|\psi_{k+1} \theta|_\R(t) - t\|_{\cC^r(I_{k+1})} <
        \e_{k+1}. \eqno(1)
$$
Apply lemma 2.3 with $\l$ replaced by $H^{-1} \l$,
$K=R\bar \B$, $r$ unchanged, $Z = T \cap I_{k+1}$, and
$\e = \min\{\d, \e_{k+1}/2\}$. This gives a shear
$\theta_{k+1}(z) = z + h_{k+1}(z_1) v_{k+1}$
with $\pi_1 v_{k+1} = 0$ such that
$$
    \cases{ |\theta_{k+1}(z)-z| < \min\{\d, \e_{k+1}/2\},\ \
         z \in\Delta_{k+1}; \cr
      \theta_{k+1}(\C) \cap H^{-1} \l(\R) = H^{-1} \l(T \cap I_{k+1}); \cr
      h_{k+1}^{(s)}(t) = 0,\ \ t \in T \cap I_{k+1},\ 0 \le s \le r,}
$$
and such that (1) holds with $\theta = \theta_{k+1}$.
Also, by the choice of $R$ and $\d$,
$$
        \psi_{k+1} \theta_{k+1}(\bar{\R \bs I_{k+1}}) \cap A
        = \emptyset.
$$
Taking $\Psi_{k+1} = \psi_{k+1} \theta_{k+1}$ and
$$
        f_{k+1} = H\circ \theta_{k+1}=
        \Phi_{k+1}\circ f_k\circ \Psi_{k+1}
$$
we obtain $(5_{k+1})$ and preserve the remaining hypotheses.
Hence we obtain $(1_{k+1})$--$(7_{k+1})$. Note that
$(k+1)\B \subset B_{k+1}$ so we also obtain (a)--(f),
thus finishing the inductive step.

\medskip

The case $k=1$ is similar to the general step.  First apply
proposition 1.2 with $K = \emptyset$, $C = [-3,3] \subset \C$, $F =
\l$, $\e = C^{-1} \inf\{\eta(t) \colon t \in [-3,3]\}$, and $Z = T \cap
[-3,3]$ to get $\phi_1^1 \in \Aut\C^n$ satisfying the conclusions of
that proposition.  Choose $R_1 \ge 2$ such that $\phi_1^1(\Delta_1)
\subset (R_1-1)\B$, choose $a_1 > 4$ to get $(6_1)$, and let $I_1 =
[-a_1, a_1]$.  Choose $\d_1$ to satisfy (f) for $j=1$.  

Define a proper $\cC^\infty$ embedding $\l_0$ as in lemma 3.1
so that (i)--(v) are satisfied with $\l_0$ in place of $\l_k$,
$\phi_1^1$ in place of $f_k$, $3$ in place of $a_k$, $[-3,3]$ in place
of $I_k$, and $\phi_1^1(\Delta_1)$ in place of $K_k$.  Apply lemma 2.4
with $\l = \l_0$, $f = \phi_1^1$, $I = [-2,2]$,
$K = \phi_1^1(\Delta_1)$, $a = a_1$, $\e = \delta_1$, and $T$ and $r$
unchanged.  This gives $\phi_1^2 \in \Aut\C^n$ such that
$$
    |\phi_1^2(z) - z| < \d_1 \le 1/2,\ \ z \in \phi_1^1(\Delta_1),
$$
and such that $\Phi_1 = \phi_1^2 \phi_1^1$ satisfies
$$
     \cases{ \|\Phi_1 - \l_0\|_{\cC^r(I_1)} < \e; \cr
     \Phi_1^{(s)}(t) = \l_0^{(s)}(t),\ \ t \in T \cap I_1,\ 0 \le s \le r.}
$$
As before, we can apply lemmas 2.2 and 2.3 to obtain $\e_1>0$ and
a shear $\Psi_1$ such that the hypotheses $(1_1)$--$(7_1)$ hold
for $f_1 = \Phi_1 \Psi_1$, and (a)--(f) hold for $k=1$. This
completes the base case.

\medskip
To finish the proof of theorem 1.1, note that
$$
        \Psi_1 \cdots \Psi_k(z) = z + \sum_{j=1}^k (g_j(z_1) e_2 +
        h_j(z_1) v_j)
$$
and that (c) implies
$$
        |g_j(z_1) e_2 + h_j(z_1) v_j| < 2^{-j},\quad |z_1|<j.
$$
Hence this sum converges uniformly on compacts to a shear
$\Psi(z) = z + G(z_1)$ for some holomorphic map
$G\colon \C \to \{0\} \times \C^{n-1}$.

By proposition 4.2 in [6], the composition
$\Phi_k \circ \cdots \circ \Phi_1$ converges locally uniformly
to a biholomorphic map from a domain $\Omega$ onto $\C^n$,
and
$$
        \Omega = \cup_{k=1}^\infty
           (\Phi_k \cdots \Phi_1)^{-1}(B_{k-1}).
$$
We claim that $\Psi(\C \times \{0\}) \subset \Omega$.  Let $k>1$.
By $(1_k)$ we have
$$
        \Psi_1 \cdots \Psi_k(\Delta_{k-1} + \e_k \bar \B) \subset
        (\Phi_k \cdots \Phi_1)^{-1}(B_{k-1}).
$$
Since $|\Psi_j(z) - z| < \e_j$ on $\Delta_{k-1}$ for $j \ge k$, and
$\sum_{j=k+1}^\infty \e_j < \e_k$ by (c), we see that
$$
        \lim_{m \to \infty} \Psi_{k+1} \cdots \Psi_m(z) \in
        \Delta_{k-1} + \e_k \bar \B
$$
for each $z \in \Delta_{k-1}$.  Hence
$$
        \Psi(\Delta_{k-1}) \subset (\Phi_k \cdots \Phi_1)^{-1}(B_{k-1})
        \subset \Omega, \qquad k>1,
$$
so the claim holds.  In particular,
$\Phi \Psi \colon \C \to \C^n$ is a \phe.

Finally, using the conditions $(1_k)$--$(7_k)$, we see that
$\Phi \Psi \colon \C \hra \C^n$ is a proper holomorphic embedding
with the desired properties.
\endpr

%
%
%
%

\bigskip
\ni\bf References. \rm
\medskip

\item{1.} Arakelian, N.\ U.:
Uniform approximation on closed sets by entire functions (Russian).
{\it Izv.\ Akad.\ Nauk SSSR} {\bf 28}, 1187--1206 (1964)

\item{2.} Buzzard, G., Forn\ae ss, J.E.:
An embedding of $\C$ into $\C^2$ with hyperbolic complement.
{\it Math.\ Ann.}, to appear

\item{3.} Carleman, T.: Sur un theorem de Weierstrass.
{\it Arkiv Mat., Astr.\ och Fysik} {\bf 20B}, no.\ 4 (1927)

\item{4.} Eliashberg, Y., Gromov, M.:
Embeddings of Stein manifolds of dimension $n$ into the affine space
of dimension $3n/2+1$.
{\it Ann.\ of Math.} {\bf 136}, 123--135 (1992)

\item{5.} Forstneric, F.: Approximation by automorphisms
on smooth submanifolds of $\C^n$.
{\it Math.\ Ann.} {\bf 300}, 719--738 (1994)

\item{6.} Forstneric, F.:
Interpolation by holomorphic automorphisms and embeddings in $C^n$.
Preprint, 1996

\item{7.} Forstneric, F., Globevnik, J., Rosay, J.-P.:
Non straightenable complex lines in $\C^2$.
{\it Arkiv Mat.} {\bf 34} (1996)
           
\item{8.} Forstneric, F., Globevnik, J., Stens\o nes, B.:
Embedding holomorphic discs through discrete sets.
{\it Math.\ Ann.} {\bf 304}, (1995)

\item{9.} Forstneric, F., L\o w, E.:
Global holomorphic equivalence of certain smooth submanifolds in $\C^n$.
Preprint, 1996

\item{10.} Gaier, D.: \it Lectures on Carleman Approximation. \rm
Birkh\"auser: Boston 1987

\item{11.} H\"ormander, L.:
\it An Introduction to Complex Analysis in Several Variables,
3rd ed. \rm North Holland: Amsterdam 1990

\item{12.}  Narasimhan, R.: {\it Analysis on real and complex
manifolds}, 3rd printing. North-Hol\-land: Amsterdam 1985

\item{13.} Rosay, J.-P., Rudin, W.:
Arakelian's approximation theorem.
{\it Amer.\ Math.\ Monthly} {\bf 96}, 432--434 (1989)

\item{14.} Rudin, W.: \it Real and Complex Analysis. \rm
McGraw Hill: New York 1970

\bigskip\rm
\noindent{\ninecsc
Gregery T. Buzzard,
Department of Mathematics, Indiana University,
Bloomington, IN 47405, USA} 
\smallskip
\noindent{\ninecsc
Franc Forstneric,
Department of Mathematics, University of Wisconsin, Madison,
WI 53706, USA}
\bigskip

\bye